# MARKET FREE LUNCH AND LARGE FINANCIAL MARKETS

By Irene Klein

*University of Vienna*

The main result of the paper is a version of the fundamental theorem of asset pricing (FTAP) for large financial markets based on an asymptotic concept of no market free lunch for monotone concave preferences. The proof uses methods from the theory of Orlicz spaces. Moreover, various notions of no asymptotic arbitrage are characterized in terms of no asymptotic market free lunch; the difference lies in the set of utilities. In particular, it is shown directly that no asymptotic market free lunch with respect to monotone concave utilities is equivalent to no asymptotic free lunch. In principle, the paper can be seen as the large financial market analogue of [*Math. Finance* **14** (2004) 351–357] and [*Math. Finance* **16** (2006) 583–588].

**1. Introduction.** In a semimartingale model of a securities market, Frittelli [5] introduced the notion of *no market free lunch* (NMFL) which is a condition of no arbitrage type induced by the preferences of all agents in the market. Preferences are given by expected utility, where the expectation is taken with respect to an equivalent probability measure which, together with the utility function, determines an agent. Frittelli gave characterizations of the well-known notions of no arbitrage (NA) and no free lunch with vanishing risk (NFLVR) in terms of NMFL, the difference between the characterizations depending on the class of utility functions. Moreover, he provided a new version of the fundamental theorem of asset pricing (FTAP), that is, he gave a direct proof of the equivalence between NMFL with respect to monotone concave utility functions and the existence of a (local/sigma) martingale measure. The proof is not based on any other version of the FTAP.

A large financial market is a sequence of traditional market models, each based on a finite number of assets. Under the assumption that there does

Received July 2005; revised March 2006.

not exist any kind of arbitrage on any of the finite markets, Kabanov and Kramkov [10] introduced the notions of no asymptotic arbitrage of first kind (NAA1) and second kind (NAA2). These notions give one-sided versions of the FTAP for large financial markets in the sense that they are equivalent to the existence of a sequence of equivalent martingale measures with a contiguity property with respect to the sequence of the original probability measures in *one* direction; see [10, 11, 16, 17].

However, the direct analogue of the existence of an equivalent martingale measure for the case of a large financial market seems to be the existence of a sequence of martingale measures which is contiguous with respect to the sequence of the original probabilities and vice versa (i.e., *bicontiguous*). So, the question naturally arises as to whether the symmetric property of bicontiguity has an economic interpretation. The general answer to this question was given in [12] where the notion of no asymptotic free lunch (NAFL) was introduced, this notion being the large financial market analogue of the concept of no free lunch (NFL) of Kreps [18]. This lead to a general version of the FTAP for large financial markets. However, the condition NAFL is rather technical and, unfortunately, cannot be replaced by a weaker condition of this type in general. For *continuous* processes, the bicontiguity property is equivalent to no asymptotic free lunch with bounded risk (NAFLBR) which has an intuitive interpretation; see [13].

In the present paper, I introduce the notion of no asymptotic market free lunch (NAMFL) which is the translation of the concept of no market free lunch of Frittelli to large financial markets. It turns out that NAMFL with respect to monotone concave preferences is also equivalent to the existence of a bicontiguous sequence of martingale measures, which is the main result of this paper. This is another version of the FTAP for large financial markets in full generality with the additional property that the notion of NAMFL has a clear economic interpretation via preferences of investors. The proof is based on some ideas of [5], but is more involved than the proof of the analogous result there. One has to introduce a uniform concept of no market free lunch with respect to monotone concave preferences (which turns out to be equivalent to the nonuniform concept) and use the theory of Orlicz spaces and a quantitative version of the Halmos–Savage theorem. The proof illustrates the close ties between the theory of monotone concave utilities and that of Orlicz spaces. This connection obviously also holds in the small market case: with the use of Orlicz spaces, one can directly prove that NMFL with respect to monotone concave preferences is equivalent to the concept of NFL. This was done in [15].

Besides the FTAP result, various no asymptotic arbitrage conditions from the literature are characterized in terms of no asymptotic market free lunch. The difference lies in the choice of the set of utilities defining the market free lunch. By the FTAP result, one may already suspect that the characterization of NAFL in terms of no market free lunch is its equivalence to NAMFL with respect to monotone concave preferences. In fact, this result is the step justifying the fact that we can replace the NAMFL condition by a uniform version. Moreover, this is an analogous result to [15] for large financial markets.

When speaking about large financial markets, one must remark on the history of the theory. The idea of asymptotic arbitrage first appeared in [22] and in the important note of Huberman [8]. The theory of large financial markets can be seen as a "modern" version of arbitrage pricing theory (APT). However, the classic APT does not touch the problem of the existence of an equivalent martingale measure, so it does not provide a version of the FTAP for a large financial market.

We conclude the Introduction with a short overview of how the paper is organized. In Section 2 we introduce some notation and deal with the small market case, that is, we compile some results on market free lunch for a traditional market model, following [5] and [15]. In Section 3 we introduce large financial markets and define an asymptotic concept of no market free lunch. In Section 4 a version of the FTAP for large financial markets based on no asymptotic market free lunch is stated and proved. In Section 5 it is shown that NAMFL with respect to monotone concave preferences does indeed hold uniformly and, moreover, is equivalent to NAFL. Theorem 4.2 and Corollary 5.3 are the main results of the paper, and together, they give a general formulation of the FTAP in terms of no market free lunch with respect to monotone concave preferences in the case of a large financial market. In Section 6 characterizations of various no asymptotic arbitrage properties are given in terms of no asymptotic market free lunch. Finally, in the Appendix, some properties of Orlicz spaces, Young functions and the Mackey topology are compiled, as well as some facts on contiguity of sequences of probabilities, a convex duality result of Bellini and Frittelli [1] based on [21] and a quantitative version of the Halmos–Savage theorem; see [17].

**2. Notation and no market free lunch on a small market.** Let $(\Omega, \mathcal{F}, (\mathcal{F}_t)_{t \in [0, \infty)}, P)$ be a filtered probability space where the filtration satisfies the usual assumptions. Let $\mathbb{P}$ be the class of probability measures on $(\Omega, \mathcal{F})$ equivalent to $P$. Whenever it is clear which $P$ is meant, the notation $L^p$ is used for $L^p(\Omega, \mathcal{F}, P)$ (where $p = 0$, $p = 1$ or $p = \infty$); whenever the dependence on a certain measure $R$ is stressed, we use the notation $L^p(R)$. Let $(S_t)_{t \in [0, \infty)}$ be an $(\mathcal{F}_t)$-adapted semimartingale with values in $\mathbb{R}^d$, describing the price processes of $d$ tradeable assets. Let $H$ be a predictable $S$-integrable process and $(H \cdot S)_t$ the stochastic integral of $H$ with respect to $S$. The process $H$ is an *admissible trading strategy* if there exist some $a > 0$ such that $(H \cdot S) \geq -a$. Define

$$\mathbf{K} = \left\{ (H \cdot S)_\infty \colon H \text{ admissible and } (H \cdot S)_\infty = \lim_{t \to \infty} (H \cdot S)_t \text{ exists a.s.} \right\},$$

$$\mathbf{C} = (\mathbf{K} - L_+^0) \cap L^\infty.$$

$\mathbf{K}$ can be interpreted as the cone of all replicable claims, and $\mathbf{C}$ as the cone of all claims in $L^\infty$ that can be superreplicated. That is, claims in $\mathbf{C}$ are dominated by elements of $\mathbf{K}$. As a piece of notation, define

$$\mathbf{K}_1 = \{ f \in \mathbf{K} \colon f \geq -1 \}.$$

Define the set of separating measures

$$M = \{ Q \ll P \colon \mathbf{K} \subseteq L^1(Q) \text{ and } E_Q[f] \leq 0 \text{ for all } f \in \mathbf{K} \}.$$

If $S$ is bounded (resp. locally bounded), then $M$ consists of all $P$-absolutely continuous probability measures such that $S$ is a martingale (resp. local martingale). In general, for unbounded $S$, $M$ is the set of $P$-absolutely continuous probabilities such that the admissible stochastic integrals are supermartingales.

REMARK 2.1. Observe that in the definition of $\mathbf{K}$, one has to be careful that $(H \cdot S)_\infty$ is well defined. But this is not a problem as later on there will always be some $Q \in M \cap \mathbb{P}$; see (3.1). Therefore, for admissible $H$, $H \cdot S$ is a $Q$-supermartingale and $(H \cdot S)_\infty$ exists a.s.

Recall the well-known conditions of no arbitrage (NA), no free lunch with vanishing risk (NFLVR) and no free lunch (NFL):

(NA)                              $\mathbf{C} \cap L_+^\infty = \{0\},$

(NFLVR)                          $\bar{\mathbf{C}} \cap L_+^\infty = \{0\},$

(NFL)                            $\bar{\mathbf{C}}^* \cap L_+^\infty = \{0\},$

where $\bar{\mathbf{C}}$ is the $L^\infty$-norm closure of $\mathbf{C}$ and $\bar{\mathbf{C}}^*$ the closure of $\mathbf{C}$ with respect to the weak-star topology.

The fundamental theorem of asset pricing (FTAP) says that an appropriate "no-arbitrage" condition is equivalent to $M \cap \mathbb{P} \neq \varnothing$. If $S$ is (locally) bounded, this means that there exists an equivalent (local) martingale measure. In the unbounded case, this implies that the set $M^\sigma$ of equivalent sigma-martingale measures is nonempty and that $M^\sigma$ is dense in $M$ for the variation topology; see [4]. For finite $\Omega$ [6], and for general $\Omega$ and a finite number of trading dates [2], NA is sufficient. For general $\Omega$ and continuous time, the FTAP was proved in full generality by Delbaen and Schachermayer [3, 4], using the stronger condition NFLVR. Moreover, it is clear that the even stronger condition NFL of Kreps [18] is also equivalent to $M \cap \mathbb{P} \neq \varnothing$.

*No market free lunch.* Frittelli [5] introduced a notion of market free lunch that depends on the preferences of the investors in the market. Let $\mathbb{U}$ be a class of utility functions $u : \mathbb{R} \to \mathbb{R} \cup \{-\infty\}$ where the domain of $u$ is $D(u) = \{x \in \mathbb{R} : u(x) > -\infty\}$. The preference $\succeq$ of a certain agent (determined by some $R \in \mathbb{P}$ and $u \in \mathbb{U}$) can be represented by expected utility:

$$f_1 \succeq f_2 \quad \Longleftrightarrow \quad E_R[u(f_1)] \geq E_R[u(f_2)].$$

A *market free lunch* with respect to $\mathbb{U}$ is a $w \in L^\infty_{++}$ such that for all $R \in \mathbb{P}$ and all $u \in \mathbb{U}$, $\sup_{f \in \mathbf{C}} E_R[u(f - w)] \geq u(0)$, where $L^\infty_{++} = \{w \in L^\infty : P(w \geq 0) = 1 \text{ and } P(w > 0) > 0\}$.

This means that all agents $(R, u)$ in the market regard the same claim $w$ as a free lunch. The supremum is computed on $\{f \in \mathbf{C} : E_R[u^-(f - w)] < +\infty\}$ and if this set is empty, the supremum is $-\infty$.

DEFINITION 2.2. There is *no market free lunch* with respect to $\mathbb{U}$, denoted by NMFL($\mathbb{U}$), if for each $w \in L^\infty_{++}$, there exist $R \in \mathbb{P}$ and $u \in \mathbb{U}$ such that $\sup_{f \in \mathbf{C}} E_R[u(f - w)] < u(0)$.

The notion of NMFL($\mathbb{U}$) depends on the class $\mathbb{U}$. Frittelli used the following sets of utility functions:

DEFINITION 2.3.

$\mathbb{U}_0 = \{u : \mathbb{R} \to \mathbb{R} \cup \{-\infty\} \text{ such that } u \text{ is nondecreasing on } \mathbb{R}\}$;

$\mathbb{U}_1 = \{u \in \mathbb{U}_0 : u \text{ is left continuous in } 0 \in \text{interior}(D(u))\}$;

$\mathbb{U}_2 = \{u \in \mathbb{U}_0 : u \text{ is finite-valued and concave on } \mathbb{R}\}$.

It is clear that NMFL($\mathbb{U}_2$) $\Rightarrow$ NMFL($\mathbb{U}_1$) $\Rightarrow$ NMFL($\mathbb{U}_0$) as $\mathbb{U}_2 \subset \mathbb{U}_1 \subset \mathbb{U}_0$. Frittelli gave characterizations of (NA) and (NFLVR) in terms of no market free lunch and proved a version of the FTAP via NMFL($\mathbb{U}_2$).

THEOREM 2.4 ([5]).

(a) *NMFL($\mathbb{U}_0$)* $\Leftrightarrow$ *NA.*
(b) *NMFL($\mathbb{U}_1$)* $\Leftrightarrow$ *NFLVR.*
(c) *NMFL($\mathbb{U}_2$)* $\Leftrightarrow$ $M \cap \mathbb{P} \neq \varnothing$ *(FTAP).*

Delbaen and Schachermayer proved that NFLVR $\Leftrightarrow M \cap \mathbb{P} \neq \varnothing$; see [3, 4]. So it is clear that, moreover, NMFL($\mathbb{U}_1$) $\Leftrightarrow$ NMFL($\mathbb{U}_2$). It will become clear later on that this is not the case for the analogous notions in large financial markets. Indeed, we will see in Section 6 that NAMFL($\mathbb{U}_1$) $\Leftrightarrow$ no asymptotic arbitrage of first kind and in Sections 4 and 5 that NAMFL($\mathbb{U}_2$) $\Leftrightarrow$ there is a sequence of measures $Q^n \in M^n \cap \mathbb{P}^n$ that is bicontiguous with respect to $(P^n)$. From the literature, we know that then NAMFL($\mathbb{U}_1$) cannot be equivalent to NAMFL($\mathbb{U}_2$); see [10, 11, 12, 16, 17].

In [15], it was shown directly that NMFL($\mathbb{U}_2$) is equivalent to NFL:

PROPOSITION 2.5.   *NMFL($\mathbb{U}_2$) $\Leftrightarrow$ NFL.*

In Section 5, a result analogous to Proposition 2.5 will be proved in the case of large financial markets; see Proposition 5.2.

**3. Large financial markets and no market free lunch.** In a *large financial market*, a sequence of small market models is considered, that is, a sequence of semimartingales $(S^n)$ where $S^n$ is based on $(\Omega^n, \mathcal{F}^n, (\mathcal{F}_t^n), P^n)$. The interpretation of the superscript $n$ in expressions such as $\mathbf{K}^n$, $\mathbf{C}^n$, $M^n$, $\mathbb{P}^n$, and so on, is then obvious. Throughout this paper, the assumption

$$(3.1) \qquad M^n \cap \mathbb{P}^n \neq \varnothing \qquad \text{for all } n \in \mathbb{N},$$

is made. This ensures that on all the small markets, NA, NFLVR and NFL hold and, moreover, NMFL($\mathbb{U}$) for $\mathbb{U} = \mathbb{U}_0, \mathbb{U}_1, \mathbb{U}_2$.

DEFINITION 3.1.   A sequence of probability measures $(R^n)$ is contiguous with respect to $(P^n)$, denoted by $(R^n) \lhd (P^n)$, if for all $A^n \in \mathcal{F}^n$, $P^n(A^n) \to 0$ implies $R^n(A^n) \to 0$.

The following notation is used: $(R^n) \lhd \rhd (P^n)$ means that the sequence of probability measures $(R^n)$ is contiguous with respect to the sequence of probability measures $(P^n)$ and vice versa.

Contiguity is used as a concept of absolute continuity of probability measures in a uniform way for sequences of probability measures. For a precise statement, see Lemma A.10 of the Appendix.

Though any form of no arbitrage holds on any of the small markets, there is still the possibility of achieving various approximations of an arbitrage profit by trading on the sequence of small markets; see Sections 5 and 6 and [10, 11, 12, 13, 16, 17].

For the definition of an asymptotic concept of no market free lunch for the large financial market which is inspired by Definition 2.2, we need the sets

$$(3.2) \qquad D^{\varepsilon, n} = \{w \in L^\infty(P^n) : 0 \leq w \leq 1 \text{ and } E_{P^n}[w] \geq \varepsilon\}.$$

DEFINITION 3.2. There is no asymptotic market free lunch for $\mathbb{U}$, denoted by NAMFL($\mathbb{U}$), if for every $\varepsilon > 0$ and every net $\mathcal{W} = (w^\gamma)_{\gamma \in \Gamma}$ with $w^\gamma \in D^{\varepsilon, n_\gamma}$, there exist $u \in \mathbb{U}$ and a sequence $R^n \in \mathbb{P}^n$ with $(R^n) \triangleleft \triangleright (P^n)$ such that for all sequences $(w^k) \subseteq \mathcal{W}$:

$$\limsup_{n \to \infty} \sup_{f^k \in \mathbf{C}^{n_k}} E_{R^{n_k}}[u(f^k - w^k)] < u(0).$$

So, if there is an asymptotic market free lunch, then all agents [represented by a utility $u \in \mathbb{U}$ and a sequence of measures $(R^n)$] regard the same net $w^n \in D^{\alpha, n}$ (for a fixed $\alpha > 0$) as a kind of asymptotic free lunch. Indeed, "in limit" they all are able to hedge the claim $\mathcal{W}$ according to their beliefs and preferences (at least on a subsequence of markets), as for all $R^n \in \mathbb{P}^n$ with $(R^n) \triangleleft \triangleright (P^n)$ and $u \in \mathbb{U}$ there is a subsequence such that $\limsup_{n \to \infty} \sup_{f^{n_k} \in \mathbf{C}^{n_k}} E_{R^{n_k}}[u(f^k - w^k)] \geq u(0)$. The $\alpha$ is fixed, so the profit does not disappear in the limit.

REMARK 3.3. For a sequence $(w^n)$ such that $w^n \in D^{\alpha, n}$ for all $n \in \mathbb{N}$, the following holds: if $(P^n) \triangleleft (R^n)$, then there exists $\beta > 0$ such that $w^n \in D^\beta(R^n)$ for all $n \in \mathbb{N}$.

Indeed, it is clear that $P^n(w^n \geq \frac{\alpha}{2}) \geq \frac{\alpha}{2}$ for each $w^n \in D^{\alpha, n}$. By Lemma A.10, there exists $\delta > 0$ such that $R^n(w^n \geq \frac{\alpha}{2}) \geq \delta$ for all $n$, and therefore $E_{R^n}[w^n] \geq \frac{\alpha}{2}\delta =: \beta$. This shows that $w^n \in D^\beta(R^n)$.

Note that in the large financial market case, the sequence of sets $D^{\varepsilon, n}$ plays the role of $L^\infty_{++}$ of Definition 2.2. This ensures that the strictly positive part of the functions does not disappear when $n \to \infty$.

## 4. A FTAP for large financial markets based on no asymptotic market free lunch.

A FTAP for a large financial market shows the equivalence result of a notion of no asymptotic arbitrage and the existence of an "equivalent martingale measure" [which is a sequence of measures $(Q^n)$ such that (i) $Q^n \in M^n \cap \mathbb{P}^n$ and (ii) $(Q^n)$ is *uniformly equivalent* with respect to the sequence $(P^n)$, i.e., $(P^n) \triangleleft \triangleright (Q^n)$]. In the case of (locally) bounded $S^n$, this means that there exists a sequence $(Q^n)$ of equivalent (local) martingale measures such that $(Q^n) \triangleleft \triangleright (P^n)$. In the unbounded case, it implies the existence of a sequence $(Q^n)$ of equivalent sigma-martingale measures such that $(Q^n) \triangleleft \triangleright (P^n)$ (this is an easy consequence of the fact that $M^n_\sigma$ is dense in $M^n$ for the variation topology, for all $n$).

In the case of continuous processes, the economically intuitive notion of no asymptotic free lunch with bounded risk (NAFLBR)—for the definition see Section 6—is necessary and sufficient; see [13]. However, in the general (i.e., noncontinuous) case, NAFLBR is not sufficient. One needs the stronger concept of no asymptotic free lunch (NAFL) which was introduced in [12]. This is the generalization of the concept of NFL by Kreps to large financial markets. The definition of NAFL is recalled in Section 5. NAFL provides a version of the FTAP for large financial markets in full generality. Unfortunately, the concept of NAFL is of a highly technical and unintuitive nature, involving subtle descriptions of Mackey neighborhoods of 0. Therefore, it is of interest to know if there is a notion that gives the same general result but which can be interpreted in an easier way. It turns out that one can prove a theorem in full generality using NAMFL*($\mathbb{U}_2$), a uniform version of NAMFL($\mathbb{U}_2$); see Theorem 4.2 below. This is the large financial market analogue of Frittelli's result [Theorem 2.4(c)]. Moreover, NAMFL*($\mathbb{U}_2$) is fairly intuitive, economically reasonable and less technical than NAFL. But we must still check whether it is admissible to use the *uniform* version of no asymptotic market free lunch. In Section 5, the gap will be closed; it will be shown there that the use of the uniform notion is indeed justified, as NAMFL*($\mathbb{U}_2$) is equivalent to NAMFL($\mathbb{U}_2$). Moreover, this gives an analogous result to Proposition 2.5 for large financial markets as it turns out that both notions are equivalent to no asymptotic free lunch (NAFL). So, in fact, Theorem 4.2 and Corollary 5.3 together give a complete formulation of the FTAP in terms of NAMFL with respect to monotone concave utilities.

In the following, some facts about Orlicz spaces $L_F(P)$ and Young functions are used. The set of all Young functions is denoted by $\mathcal{Y}$. All relevant definitions and results are compiled in the Appendix.

DEFINITION 4.1. There is *uniformly* no asymptotic market free lunch with respect to $\mathbb{U}_2$ on the large financial market, NAMFL*($\mathbb{U}_2$) if for every $\varepsilon > 0$, there is $\delta > 0$, $u \in \mathbb{U}_2$ and a sequence $R^n \in \mathbb{P}^n$ with $(R^n) \lhd \rhd (P^n)$ such that for every sequence $w^n \in D^{\varepsilon,n}$,

$$\limsup_{n \to \infty} \sup_{f^n \in \mathbf{C}^n} E_{R^n}[u(f^n - w^n)] < u(0) - \delta.$$

This is called *uniformly* since the $u$ and the $(R^n)$ in the NAMFL′ condition depend only on $\varepsilon$ and not on the special choice of a sequence $w^n \in D^{\varepsilon,n}$. On the other hand, in the NAMFL condition for each net $(w^\gamma)_{\gamma \in \Gamma}$ with $w^\gamma \in D^{\varepsilon,n_\gamma}$, there is some $u \in \mathbb{U}_2$ and some sequence $(R^n)$. It will become clear that, in contrast to the small market case, the uniformity is crucial when one is dealing with large financial markets. In the proof, we will have to use the fact that the $u$ and $(R^n)$ only depend on $\varepsilon$. Without this uniformity, the proof of the ($\Rightarrow$) part of Theorem 4.2 would be very hard, or even not possible.

THEOREM 4.2 (A version of the FTAP for large financial markets).

*NAMFL'*($\mathbb{U}_2$) $\Leftrightarrow$ *there exists a sequence of measures* $Q^n \in M^n \cap \mathbb{P}^n$ *such that* $(Q^n) \triangleleft \triangleright (P^n)$.

PROOF. ($\Leftarrow$) Let $(Q^n) \triangleleft \triangleright (P^n)$, fix $\varepsilon > 0$ and take $u(x) = x$. Take an arbitrary $w^n \in D^{\varepsilon,n}$. Then $E_{P^n}[w^n] \geq \varepsilon$ and, as $w^n \in D^{\varepsilon,n}$ and $(P^n) \triangleleft (Q^n)$, $E_{Q^n}[w^n] > \delta$ for some $\delta > 0$ that depends only on $\varepsilon$ (see Remark 3.3). So, as $Q^n \in M^n$, $E_{Q^n}[f^n] \leq 0$ and

$$\sup_{f^n \in \mathbf{C}^n} E_{Q^n}[u(f^n - w^n)] = \sup_{f^n \in \mathbf{C}^n} E_{Q^n}[f^n] - E_{Q^n}[w^n] < -\delta,$$

hence NAMFL*($\mathbb{U}_2$). ($\Rightarrow$) Let $\varepsilon > 0$ be fixed. By NAMFL*($\mathbb{U}_2$), there exist $\delta > 0$, $\tilde{u} \in \mathbb{U}_2$ and a sequence $R^n \in \mathbb{P}^n$ with $(R^n) \triangleleft \triangleright (P^n)$ such that for every sequence $\mathbf{A} = (A^n)$ with $A^n \in \mathcal{F}^n$ and $P^n(A^n) \geq \varepsilon$ (i.e., $\mathbb{1}_{A^n} \in D^{\varepsilon,n}$),

$$\limsup_{n \to \infty} \sup_{f^n \in \mathbf{C}^n} E_{R^n}[\tilde{u}(f^n - \mathbb{1}_{A^n})] < \tilde{u}(0) - 3\delta.$$

Without loss of generality, let $\tilde{u}(0) = 0$ [just replace $\tilde{u}(x)$ by $\tilde{u}(x) - \tilde{u}(0)$ and the same holds]. By Lemma A.4, for $\delta > 0$ as above, there exist $u \in \mathbb{U}_2$ such that $u(x) \leq \tilde{u}(x) + \delta$ and there exist a Young function $F$ such that $u$ is defined as $u_F$ in (A.2). It is clear that for this $u$,

$$(4.1) \quad \limsup_{n \to \infty} \sup_{f^n \in \mathbf{C}^n} E_{R^n}[u(f^n - \mathbb{1}_{A^n})] < u(0) - 2\delta = -2\delta \leq 0 = u(+\infty).$$

So the assumptions of Theorem A.14 are fulfilled for every $n \geq N^{\varepsilon, \mathbf{A}}$ for $G = \mathbf{C}^n$, $N = M^n$, $w = -\mathbb{1}_{A^n}$ and $u$. This implies that there exist $Q^{A^n} \in M^n$ and $\lambda^{A^n} > 0$ such that, for all $n \geq N^{\varepsilon, \mathbf{A}}$,

$$(4.2) \quad E_{R^n}\left[-\lambda^{A^n} \frac{dQ^{A^n}}{dR^n} \mathbb{1}_{A^n} + v\left(\lambda^{A^n} \frac{dQ^{A^n}}{dR^n}\right)\right]$$
$$= \sup_{f^n \in \mathbf{C}^n} E_{R^n}[u(f^n - \mathbb{1}_{A^n})] < -\delta,$$

where $v(y) = \sup_{x \in \mathbb{R}}(u(x) - xy)$.

CLAIM 1. *There exists some* $N < \infty$ *such that* $N^{\varepsilon, \mathbf{A}} \leq N$ *for all sequences* $\mathbf{A}$ *as above.*

Indeed, otherwise for each $k$, there is a sequence $(A_k^n)$ such that for some $n_k \geq k$,

$$(4.3) \quad \sup_{f^{n_k} \in \mathbf{C}^{n_k}} E_{R^{n_k}}[u(f^{n_k} - \mathbb{1}_{A_k^{n_k}})] \geq -\delta.$$

Define a new sequence $B^k = A_k^{n_k}$. Then, still, $P^{n_k}(B^k) \geq \varepsilon$. For the sequence $(B^k)$, (4.3) gives a contradiction to (4.2), verifying Claim 1.

From now on, fix the $N$ as above.

Claim 2. *There exist $\lambda_0 > 0$ and $\lambda_1 > 0$ such that for all sequences* **A** *as above, $\lambda_0 \leq \lambda^{A^n} \leq \lambda_1$, for all $n \geq N$. In particular, $\lambda_0 = \delta$.*

By definition, clearly $v(y) \geq u(0) = 0$ for all $y$. Hence, as $Q^{A^n}(A^n) \leq 1$ and by (4.2),

$$-\lambda^{A^n} \leq -\lambda^{A^n} Q^{A^n}(A^n) + E_{R^n}\left[v\left(\lambda^{A^n}\frac{dQ^{A^n}}{dR^n}\right)\right] < -\delta,$$

and therefore $\lambda^{A^n} > \delta =: \lambda_0$.

As for the second assertion, note that the definition of $u$ implies that $v$ is a Young function (Lemma A.5). By Lemma A.3, $\lim_{y\to\infty}\frac{v(y)}{y} = +\infty$ and therefore, a fortiori, $\lim_{y\to\infty}(v(y)-y) = +\infty$. As $0 \leq Q^{A^n}(A^n) \leq 1$, by Jensen's inequality and (4.2),

$$-\lambda^{A^n} + v(\lambda^{A^n}) = -\lambda^{A^n} + v\left(E_{R^n}\left[\lambda^{A^n}\frac{dQ^{A^n}}{dR^n}\right]\right)$$

$$\leq -\lambda^{A^n}Q^{A^n}(A^n) + E_{R^n}\left[v\left(\lambda^{A^n}\frac{dQ^{A^n}}{dR^n}\right)\right] < -\delta.$$

As $v(\lambda^{A^n}) - \lambda^{A^n} < -\delta$, there exists $\lambda_1 > 0$ (depending only on $v$ and therefore on $\varepsilon$, but not on the sequence **A**) such that $\lambda^{A^n} \leq \lambda_1$, proving Claim 2.

Thus far it has been shown that

$$(4.4) \qquad -\lambda^{A^n}Q^{A^n}(A^n) + E_{R^n}\left[v\left(\lambda^{A^n}\frac{dQ^{A^n}}{dR^n}\right)\right] < -\delta$$

for all $n \geq N$, and that $\delta \leq \lambda^{A^n} \leq \lambda_1$, where $v(y)$, $\delta > 0$, $\lambda_1 > 0$ and $N$ depend only on $\varepsilon$ and not on the special choice of the sequence $A^n$ with $P^n(A^n) \geq \varepsilon$.

(4.4) implies that there exists $\gamma > 0$ (depending only on $\varepsilon$) such that $Q^{A^n}(A^n) > \gamma$ for all $n \geq N$. Indeed, as $v \geq 0$,

$$-\delta > -\lambda^{A^n}Q^{A^n}(A^n) + E_{R^n}\left[v\left(\lambda^{A^n}\frac{dQ^{A^n}}{dR^n}\right)\right] \geq -\lambda^{A^n}Q^{A^n}(A^n),$$

which implies that $Q^{A^n}(A^n) > \frac{\delta}{\lambda^{A^n}} \geq \frac{\delta}{\lambda_1} =: \gamma$.

Let $G(y) = \frac{v(\delta y)}{\lambda_1 - \delta}$. By Lemma A.5, $v \in \mathcal{Y}$ and so, clearly, $G \in \mathcal{Y}$ as well.

Claim 3. $\frac{dQ^{A^n}}{dR^n} \in B^G(R^n) = \{f \in L_G(R^n) : \|f\|_G \leq 1\}$ *for each sequence* $(A^n)$ *as above.*

$G \in \mathcal{Y}$ and is therefore strictly increasing. Hence, as $\delta \leq \lambda^{A^n}$ and by (4.4),

$$E_{R^n}\left[G\left(\frac{dQ^{A^n}}{dR^n}\right)\right] \leq \frac{1}{\lambda_1 - \delta}E_{R^n}\left[v\left(\lambda^{A^n}\frac{dQ^{A^n}}{dR^n}\right)\right]$$

$$(4.5) \qquad < \frac{-\delta + \lambda^{A^n} Q^{A^n}(A^n)}{\lambda_1 - \delta}$$

$$\leq \frac{-\delta + \lambda_1}{\lambda_1 - \delta} = 1.$$

So $\frac{dQ^{A^n}}{dR^n} \in B^G(R^n)$ for each $n \geq N$ and each sequence $(A^n)$ as above. This proves Claim 3.

Let us recapitulate what we have found. For fixed $\varepsilon > 0$, there is a sequence $(R^n)$ [with $(R^n) \triangleleft \triangleright (P^n)$ and $R^n \in \mathbb{P}^n$] and $\gamma > 0$ and $N \in \mathbb{N}$ such that for all $n \geq N$, the following holds. For all $A^n \in \mathcal{F}^n$ with $P^n(A^n) \geq \varepsilon$, there exists $Q^{A^n} \in M^n$ with:

(a) $Q^{A^n}(A^n) > \gamma$, and

(b) $Q^{A^n} \in \mathcal{P}^n$, where $\mathcal{P}^n = M^n \cap \{Q \ll R^n : \frac{dQ}{dR^n} \in B^G(R^n)\}$. Note that $\mathcal{P}^n$ is a convex set of $P^n$-absolutely continuous probability measures (as $R^n \in \mathbb{P}^n$).

Apply Proposition A.15 for each fixed $\varepsilon$ and the corresponding $\gamma$, each $n \geq N = N(\varepsilon)$ and the convex set $\mathcal{P}^n$ of $P^n$-absolutely continuous probability measures. This gives, for each $\varepsilon > 0$, a sequence $(Q^{\varepsilon,n})_{n \geq N}$ such that $Q^{\varepsilon,n}(B^n) \geq \mu := \frac{\varepsilon^2 \gamma}{2}$ for each $B^n$ with $P^n(B^n) \geq 4\varepsilon$. Moreover, $Q^{\varepsilon,n} \in \mathcal{P}^n$; note that by Lemma A.12, this implies that $(Q^{\varepsilon,n})_{n \geq N} \triangleleft (R^n)_{n \geq N}$ and hence that $(Q^{\varepsilon,n})_{n \geq N} \triangleleft (P^n)_{n \geq N}$ for each $\varepsilon$ [as $(R^n) \triangleleft \triangleright (P^n)$].

Define, for $\varepsilon = 2^{-j}$, $j = 1, 2, \ldots$, and for $n \geq N$, $Q^n = \sum_{j=1}^{\infty} 2^{-j} Q^{2^{-j}, n}$.

CLAIM 4. *The sequence $(Q^n)_{n \geq N}$ satisfies:*

1. $Q^n \in M^n \cap \mathbb{P}^n$,
2. $(Q^n)_{n \geq N} \triangleleft \triangleright (P^n)_{n \geq N}$.

The proof can be found in [12], for the convenience of the reader, it is recalled in the Appendix (Lemma A.13). Now, for each $1 \leq n \leq N - 1$, take an arbitrary $Q^n \in M^n \cap \mathbb{P}^n \neq \varnothing$. As $Q^n \sim P^n$ for $n = 1, \ldots, N-1$, these $Q^n$ can be used to complete the sequence and we still have $(Q^n)_{n \geq 1} \triangleleft \triangleright (P^n)_{n \geq 1}$. This completes the proof. $\square$

**5. A characterization of no asymptotic free lunch in terms of no market free lunch; NAMFL with respect to $\mathbb{U}_2$ holds uniformly.** Recall the definition of NAFL from [12]. Let $F$ be a Young function and define $B^F(P^n) = \{f \in L_F(P^n) : \|f\|_F \leq 1\}$; see (A.3). By Remark A.8, $B^F(P^n)$ is a weakly compact subset of $L^1(P^n)$. As in (A.4), define the polar

$$V^F(P^n) = (B^F(P^n))^{\circ}$$

$$= \{g \in L^{\infty}(P^n) : |E_{P^n}[gh]| \leq 1 \text{ for all } h \in B^F(P^n)\}.$$

The notation $B^{F,n}$ and $V^{F,n}$ is used for $B^F(P^n)$ and $V^F(P^n)$, respectively, if it is clear which measure $P^n$ is meant.

**DEFINITION 5.1.** There is no asymptotic free lunch (NAFL) if for each $\varepsilon > 0$, there is some $F \in \mathcal{Y}$ such that $\mathbf{C}^n \cap (D^{\varepsilon,n} + V^{F,n}) = \varnothing$ for all $n \in \mathbb{N}$, where $D^{\varepsilon,n}$ is defined as in (3.2).

Definition 5.1 implies that the set $\mathbf{C}^n$ is, for each $\varepsilon > 0$, separated from $D^{\varepsilon,n}$ by some Mackey neighborhood $V^{F,n}$ of 0. See Section A.3 for an explanation of why a fundamental system for the Mackey neighborhoods of 0 [in $L^\infty(P^n)$] can be given by the sets $V^{F,n}$. Therefore, NAFL is a direct translation of NFL to a sequence of probability spaces. Indeed, if NAFL holds, it is not possible to approximate a strictly positive gain by elements of the sequence of sets $(\mathbf{C}^n)$ in a Mackey (or, equivalently, weak-star) sense.

Note that in [12], the weakly compact subsets of $L^1$ were described in a slightly different way. However, the short note [14] shows that the sets $B^{F,n}$ can be used instead of the sets $K^{\varphi,n}$ for [12]. In fact, the use of Orlicz space methods makes the proof for [12] more transparent.

The analogue of Proposition 2.5 for large financial markets will now be proven. This result shows that in Section 4, it is justified to use the uniform version of the no market free lunch condition with respect to $\mathbb{U}_2$. Moreover, it gives a characterization in terms of no market free lunch of NAFL.

**PROPOSITION 5.2.** *The following conditions are equivalent:*

  (i) *NAMFL($\mathbb{U}_2$),*
  (ii) *NAMFL*($\mathbb{U}_2$),*
  (iii) *NAFL.*

It is now clear that we can formulate a version of the FTAP with the use of NAMFL($\mathbb{U}_2$) which, together with Theorem 4.2, is the main result this the paper.

**COROLLARY 5.3.** *NAMFL($\mathbb{U}_2$) $\Leftrightarrow$ there exists a sequence of measures $Q^n \in M^n \cap \mathbb{P}^n$ such that $(Q^n) \triangleleft \triangleright (P^n)$.*

The following lemma shows that NAFL does not depend on the choice of a sequence $(R^n)$ such that $R^n \in \mathbb{P}^n$ and $(R^n) \triangleleft \triangleright (P^n)$:

**LEMMA 5.4.** *Suppose that $(f^F)_{F \in \mathcal{Y}}$ with $f^F \in \mathbf{C}^{n(F)}$ is an asymptotic free lunch with respect to $(P^n)$. Then $(f^F)_{F \in \mathcal{Y}}$ is an asymptotic free lunch with respect to any $(R^n) \triangleleft \triangleright (P^n)$, where $R^n \in \mathbb{P}^n$.*

Proof. By assumption, there exists some $\alpha > 0$ such that for all $F \in \mathcal{Y}$, there exists $f^F \in C^{n(F)}$ which can be written as $f^F = w^F + g^F$, where $w^F \in D^\alpha(P^{n(F)})$ and $g^F \in V^F(P^{n(F)})$. It will be shown that $(f^F)$ is an asymptotic free lunch for every sequence of probability measures $(R^n)$ such that $(R^n) \lhd\rhd (P^n)$.

Indeed, as $(P^n) \lhd (R^n)$, there exists $\beta > 0$ such that $w^F \in D^\beta(R^n)$ (Remark 3.3). As $(R^n) \lhd (P^n)$, Lemma A.12 implies that there exists $\psi \in \mathcal{Y}$ such that, for all $n$, $\frac{dR^n}{dP^n} \in B^\psi(P^n)$. Take now an arbitrary $\tilde{F} \in \mathcal{Y}$. Let $\tilde{G}$ be the complementary Young function of $\tilde{F}$. Moreover, let $\varphi$ be the complementary Young function of $\psi$. Defining $G(y) = (\varphi \circ 2\tilde{G})(2y)$, it is clear that $G \in \mathcal{Y}$. Let $F$ be the complementary Young function of $G$. We claim that the $g^F \in V^F(P^{n(F)})$, which is given by AFL for the Young function $F$, is in $V^{\tilde{F}}(R^{n(F)})$. Setting $n(\tilde{F}) := n(F)$, this then readily implies AFL for the sequence $(R^n)$.

Indeed, let $n = n(F)$. By Lemma A.9, $V^F(P^n) \subseteq B^G(P^n) \cap L^\infty(P^n)$. As $g^F \in V^F(P^n)$,

$$E_{P^n}[\varphi(2\tilde{G}(2|g^F|))] = E_{P^n}[G(|g^F|)] \leq 1.$$

This gives $\|\tilde{G}(2|g^F|)\|_{L_\varphi(P^n)} \leq \frac{1}{2}$. Moreover, as $\frac{dR^n}{dP^n} \in B^\psi(P^n)$, Lemma A.2(iii) implies

$$E_{R^n}[\tilde{G}(2|g^F|)] = E_{P^n}\left[\frac{dR^n}{dP^n}\tilde{G}(2|g^F|)\right]$$

$$\leq 2\left\|\frac{dR^n}{dP^n}\right\|_{L_\psi(P^n)}\|\tilde{G}(2|g^F|)\|_{L_\varphi(P^n)} \leq 1.$$

But this means that $g^F \in \frac{1}{2}B^{\tilde{G}}(R^n) \cap L^\infty(R^n) \subseteq V^{\tilde{F},n}(R^n)$, again by Lemma A.9. □

Proof of Proposition 5.2. It is clear that (ii) $\Rightarrow$ (i).

(iii) $\Rightarrow$ (ii) Suppose there is an AMFL$'(\mathbb{U}_2)$. There exists $\varepsilon > 0$ such that for all $\delta > 0$, all $u \in \mathbb{U}_2$ and all $(R^n) \lhd\rhd (P^n)$ [in particular for the sequence $(P^n)$], there exists a sequence $w^n \in D^{\varepsilon,n}$ with $\limsup_{n\to\infty}\sup_{f^n \in \mathbf{C}^n} E_{P^n}[u(f^n - w^n)] \geq u(0) - \delta$.

It will be shown that there is AFL. Indeed, let $G \in \mathcal{Y}$, arbitrary and $B^G(P^n)$ be defined as usual; see (A.3). Let $F$ be the complementary Young function of $G$ and

$$u^F(x) := \begin{cases} -F(-2x), & \text{for } x \leq 0, \\ 0, & \text{for } x > 0. \end{cases}$$

By AMFL$'(\mathbb{U}_2)$ for $\delta = \frac{1}{2}$, for $u^F \in \mathbb{U}_2$ and for the sequence $(P^n)$, there exist $n = n(F,\delta)$ and $f^n \in C^n$, and $w^n \in D^{\varepsilon,n}$ such that $E_{P^n}[u^F(f^n - w^n)] \geq -1$.

Set $g^n = (f^n - w^n)^-$. Then clearly $E_{P^n}[F(2g^n)] \le 1$ and therefore $\|g^n\|_F \le \frac{1}{2}$. By Lemma A.2(iii), for all $h \in B^G(P^n)$,

$$|E_{P^n}[g^n h]| \le E_{P^n}[|g^n h|] \le 2\|g^n\|_F \|h\|_G \le 1,$$

which means that $g^n \in (B^{G,n})^\circ = V^{G,n}$. This implies that there is an asymptotic free lunch. Indeed, if we define $\tilde{f}^n = f^n - (f^n - w^n)^+ \in \mathbf{C}^n$, then $\tilde{f}^n = w^n - g^n$ and so $\tilde{f}^n \in (D^{\varepsilon,n} + V^{G,n}) \cap \mathbf{C}^n$. This can be done for every $G \in \mathcal{Y}$, which implies AFL.

(i) $\Rightarrow$ (iii) Suppose there exists an AFL, that is, there is $\alpha > 0$ such that for each $F \in \mathcal{Y}$, there is $f^F = w^F + g^F \in C^{n(F)}$, where $w^F \in D^{\alpha,n(F)}$ and $g^F \in V^{F,n(F)}$. By Lemma 5.4, this is an AFL for an arbitrary sequence $(R^n) \triangleleft \triangleright (P^n)$ [from now on, fix $(R^n)$ and assume that the sets $D^{\alpha,n}$ and $V^{F,n}$ are defined with respect to $R^n$]. Now take any $\tilde{u} \in \mathbb{U}_2$ and w.l.o.g., assume that $\tilde{u}(0) = 0$. By Lemma A.4, for all $k$, there exists $u^k \in \mathbb{U}_2$ with $u^k(x) \le \tilde{u}(x) + 2^{-k}$ and there exists $F^k \in \mathcal{Y}$ such that $u^k := u^{F^k}$, defined as in (A.2). Moreover, because of the definition of $\mathbf{C}^{n(F)}$ which allows the subtraction of positive elements of $L^\infty$, one can assume that the $g^F \in V^{F,n}$ [where $n = n(F)$] satisfy $g^F \le 0$.

It is also clear that $\tilde{F}^k(x) := kF^k(x) \in \mathcal{Y}$. Let $G^k$ be the complementary Young function of $\tilde{F}^k$. Take the $f^k := f^{G^k} = w^k + g^k$ given by the AFL, so that $g^k \in V^{G^k, n(G^k)}$. By Lemma A.9, $V^{G^k, n_k} \subseteq B^{\tilde{F}^k, n_k}$ and so

$$E_{R^{n_k}}[\tilde{F}^k(-g^k)] = E_{R^{n_k}}[kF^k(-g^k)] \le 1.$$

This gives $E_{R^{n_k}}[F^k(-g^k)] \le \frac{1}{k}$ and $E_{R^{n_k}}[u^k(g^k)] \ge -\frac{1}{k}$. For $\tilde{u} \ge u^k - 2^{-k}$,

$$\limsup_{k \to \infty} E_{R^{n_k}}[\tilde{u}(f^k - w^k)] \ge \limsup_{k \to \infty} (E_{R^{n_k}}[u^k(f^k - w^k)] - 2^{-k})$$

$$= \limsup_{k \to \infty} (E_{R^{n_k}}[u^k(g^k)] - 2^{-k})$$

$$\ge \lim_{k \to \infty} \left( -\frac{1}{k} - 2^{-k} \right) = 0 = \tilde{u}(0).$$

This gives an AMFL($\mathbb{U}_2$).  $\square$

## 6. Characterizations of various no asymptotic arbitrage conditions in terms of no market free lunch.
The results of this section are in the style of Frittelli; see Theorem 2.4(a) and (b). Let us briefly recall some concepts of asymptotic arbitrage/free lunch that appear in the literature:

DEFINITION 6.1.  On the large financial market there is:

1. an asymptotic arbitrage of first kind (AA1) if there exist $c_k > 0$ with $c_k \to 0$ and $\xi^k \in c_k \mathbf{K}_1^{n_k}$ such that $\lim_{k \to \infty} P^{n_k}(\xi^k \ge L_k) > 0$, where $L_k > 0$ tends to infinity,

2. an asymptotic arbitrage of second kind (AA2) if there exist $\xi^k \in \mathbf{K}_1^{n_k}$ and $\alpha > 0$ such that $\lim_{k \to \infty} P^{n_k}(\xi^k \geq \alpha) = 1$,

3. a strong asymptotic arbitrage (SAA) if there exist $c_k > 0$ with $c_k \to 0$ and $\xi^k \in c_k \mathbf{K}_1^{n_k}$ and $\alpha > 0$, such that $\lim_{n \to \infty} P^{n_k}(\xi^k \geq \alpha) = 1$,

4. an asymptotic free lunch with bounded risk (AFLBR) if there exist $\alpha > 0$ and $\xi^k \in \mathbf{K}_1^{n_k}$ such that:

   (i) $P^{n_k}(\xi^k \geq \alpha) \geq \alpha$ for all $k \in \mathbb{N}$ and
   (ii) $\lim_{k \to \infty} P^{n_k}(\xi^k < -\varepsilon) = 0$ for all $\varepsilon > 0$.

A large financial market satisfies no asymptotic arbitrage of first kind (NAA1), of second kind (NAA2), no strong asymptotic arbitrage (NSAA), no asymptotic free lunch with bounded risk (NAFLBR) if there does not exist a sequence as in parts 1–4, respectively.

There are various results relating the notions NAA1 and NAA2 to the existence of a sequence of martingale measures $(Q^n)$ which is contiguous with respect to $(P^n)$ in one direction. Indeed, NAA1 $\Leftrightarrow$ $(P^n)$ is contiguous with respect to the sequence of upper envelopes of the measures $Q^n \in M^n \cap \mathbb{P}^n$, whereas NAA2 $\Leftrightarrow$ the sequence of lower envelopes of the measures $Q^n \in M^n \cap \mathbb{P}^n$ is contiguous with respect to $(P^n)$; see [11] and, in a different formulation, [16] and [17]. These results can be considered *one-sided* versions of the FTAP. As already mentioned in Section 4, in the case of continuous processes, NAFLBR is necessary and sufficient for the existence of a sequence of measures $(Q^n)$ with $Q^n \in M^n \cap \mathbb{P}^n$ and $(Q^n) \triangleleft \triangleright (P^n)$. In general, NAFLBR is not sufficient for this condition—one needs the stronger notion of NAFL; see Section 5.

REMARK 6.2. If there exist a sequence $\xi^k \in \mathbf{K}^{n_k}$ realizing AA1, AA2, SAA, AFLBR, respectively, then there exist a sequence $f^k \in \mathbf{C}^{n_k}$ with the respective properties as $\xi^k$.

Indeed, $g^{k,l} = \xi^k \vee l \to \xi^k$ in probability and $g^{k,l} = \xi^k - (\xi^k - g^{k,l}) \in \mathbf{C}^{n_k}$.

Some more sets of utility functions will now be defined; see also Definition 2.3.

DEFINITION 6.3.

$$\overline{\mathbb{U}}_0 = \{u \in \mathbb{U}_0 : 0 \in \text{interior}(D(u))\},$$

$$\mathbb{U}_3 = \{u \in \mathbb{U}_0 : u \text{ is concave and finite-valued on } [-1, \infty)\},$$

$$\mathbb{U}_4 = \{u \in \mathbb{U}_0 : u(-1) > -\infty\},$$

$$\mathbb{U}_5 = \{u : \mathbb{R} \to \mathbb{R} \cup \{-\infty\}\},$$

Note that in the list below, there is no result concerning NAMFL($\mathbb{U}_2$), as we dealt with this notion already in Section 5.

**Theorem 6.4.** *The following relations hold:*

(a) *NAMFL($\overline{\overline{\mathbb{U}}}_0$) $\Rightarrow$ NSAA,*

(b) *NAMFL($\mathbb{U}_1$) $\Leftrightarrow$ NAA1,*

(c) *NAMFL($\mathbb{U}_3$) $\Leftrightarrow$ NAFLBR,*

(d) *NAMFL($\mathbb{U}_4$) $\Rightarrow$ NAA2,*

(e) *Neither NSAA nor NAA2 imply NAMFL($\mathbb{U}_5$) which is the weakest possible form of NAMFL($\mathbb{U}$).*

**Proof.** (a) ($\Rightarrow$) Suppose there is an SAA. By definition and by Remark 6.2, this implies that there exists $f^k \in \mathbf{C}^{n_k}$ with $f^k \geq -c_k$ and there exists $\alpha > 0$ such that $\lim_{k\to\infty} P^{n_k}(A_k) = 1$, where $A_k = \{f^k \geq \alpha\}$. If we define $w^k = \alpha \mathbb{1}_{A_k}$, then $w^k \in D^{\alpha/2, n_k}$ for $k$ sufficiently large. If we let $u \in \overline{\overline{\mathbb{U}}}_0$, then there exists $\varepsilon > 0$ such that $u(x) > -\infty$ for all $x > -\varepsilon$. So, for $k$ sufficiently large (such that $c_k < \varepsilon$) it follows that, for all $(R^n) \lhd \rhd (P^n)$ with $R^n \in M^n \cap \mathbb{P}^n$,

$$\sup_{\tilde{f}^k \in \mathbf{C}^{n_k}} E_{R^{n_k}}[u(\tilde{f}^k - w^k)] \geq E_{R^{n_k}}[u(f^k - w^k)]$$

$$\geq E_{R^{n_k}}[u(f^k \mathbb{1}_{A_k^c})]$$
$$\geq u(-c_k)R^{n_k}(A_k^c) + u(0)R^{n_k}(A_k),$$

where $A_k^c = \Omega^{n_k} \setminus A_k$. Passing to the limit, this gives a contradiction to NAMFL($\overline{\overline{\mathbb{U}}}_0$) since $(R^n) \lhd \rhd (P^n)$ implies that $R^{n_k}(A_k) \to 1$.

(b) ($\Rightarrow$) If we suppose there is an AA1, then there exist $f^k \in \mathbf{C}^{n_k}$ with $f^k \geq -c_k$ for $c_k \to 0$. Moreover, there exists $\alpha > 0$ such that $P^{n_k}(A_k) \geq \alpha$ for all $k$ sufficiently large, where $A_k = \{f^{n_k} \geq \alpha\}$. If we define $w^k = \alpha \mathbb{1}_{A_k}$, then it is clear that $w^k \in D^{\alpha^2, n_k}$. Let $u \in \mathbb{U}_1$. Because of the left continuity at 0, for each $k$, there exists $\delta_k > 0$ such that $u(x) > u(0) - \frac{1}{k}$ for $-\delta_k \leq x \leq 0$. Choose a subsequence, again denoted by $k$, such that $f^k \geq w^k - \delta_k \mathbb{1}_{A_k^c}$, whence $f^k - w^k \geq -\delta_k$.

Define $\tilde{f}^k = f^k - (f^k - w^k)^+ \in \mathbf{C}^{n_k}$. Then $0 \geq \tilde{f}^k - w^k \geq -\delta_k$. Therefore, $u(\tilde{f}^k - w^k) \geq u(0) - \frac{1}{k}$. This gives a contradiction to NAMFL($\mathbb{U}_1$).

($\Leftarrow$) Suppose there is an AMFL($\mathbb{U}_1$). Define $u_k \in \mathbb{U}_1$ by

$$u_k(x) = \begin{cases} -\infty, & \text{for } x < -\dfrac{1}{k}, \\ 0, & \text{for } x \geq -\dfrac{1}{k}. \end{cases}$$

Then there exist $w^k \in D^{\alpha, n_k}$ and $f^k \in \mathbf{C}^{n_k}$ such that $E_{P^{n_k}}[u_k(f^k - w^k)] \geq -2^{-k}$. This implies that $P^{n_k}(f^k - w^k \geq -\frac{1}{k}) = 1$. For $k$ sufficiently large that $\frac{1}{k} < \frac{\alpha}{4}$, it is clear that $P^{n_k}(f^k \geq \frac{\alpha}{4}) \geq P^{n_k}(w^k \geq \frac{\alpha}{2}) \geq \frac{\alpha}{2}$, giving an AA1; take, for instance, $\tilde{f}^k = \sqrt{k} f^k \in \mathbf{C}^{n_k}$.

(c) $\Rightarrow$ If we suppose there is an AFLBR, then there exist $f^k \in \mathbf{C}^{n_k}$ such that $f^k \geq -1$. Moreover, there exists $\alpha > 0$ such that $P^{n_k}(f^k \geq \alpha) \geq \alpha$ and $\lim_{k \to \infty} P^{n_k}(f^k < -\varepsilon) = 0$ for all $\varepsilon > 0$. If we take any $u \in \mathbb{U}_3$, then there exists $\delta_n > 0$ such that $u(x) > u(0) - \frac{1}{n}$ for $-\delta_n \leq x \leq 0$. Choose $\delta_k \to 0$ and a subsequence of $f^k$, again denoted by $k$ such that $P^{n_k}(A_k) = 1 - \varepsilon_k$, where $A_k = \{f^k \geq -\delta_k\}$ and $\lim_{k \to \infty} \varepsilon_k = 0$. Let $w^k = \alpha \mathbb{1}_{\{f^k \geq \alpha\}}$ and define $\tilde{f}^k = f^k - (f^k - w^k)^+ \in \mathbf{C}^{n_k}$. Then, for every $(R^n) \triangleleft \triangleright (P^n)$,

$$E_{R^{n_k}}[u(\tilde{f}^k - w^k)] = E_{R^{n_k}}[u(-(f^k - w^k)^-)]$$
$$\geq \left( u(0) - \frac{1}{k} \right) R^{n_k}(A_k) + u(-1) R^{n_k}(A_k^c).$$

Passing to the limit for $k \to \infty$, this is a contradiction to NAMFL($\mathbb{U}_3$) since $R^{n_k}(A_k^c) \to 0$ for every $(R^n) \triangleleft \triangleright (P^n)$.

($\Leftarrow$) Suppose now there is an AMFL($\mathbb{U}_3$). Define $u_k \in \mathbb{U}_3$ by

$$u_k(x) = \begin{cases} -\infty, & \text{for } x < -1, \\ -e^{-k^2 x} + 1, & \text{for } -1 \leq x < 0, \\ 0, & \text{for } x \geq 0. \end{cases}$$

By assumption, there are $w^k \in D^{\alpha, n_k}$ and $f^k \in \mathbf{C}^{n_k}$ with $E_{P^{n_k}}[u_k(f^k - w^k)] \geq -2^{-k}$. This implies that $P^{n_k}(f^k \geq -1) = 1$ for all $k$. We now claim that $\lim_{k \to \infty} P^{n_k}(A_k) = 0$, where $A_k = \{f^k - w^k < -\frac{1}{k}\}$. Suppose to the contrary, that there exists $\beta > 0$ such that for a subsequence we have that $P^{n_k}(A_k) \geq \beta$ for all $k$. Then

$$E_{P^{n_k}}[u_k(f^k - w^k)]$$
$$= E_{P^{n_k}}[u_k(f^k - w^k) \mathbb{1}_{A_k}] + E_{P^{n_k}}[u_k(f^k - w^k) \mathbb{1}_{A_k^c}] < (-e^k + 1)\beta,$$

which is a contradiction for sufficiently large $k$. But this gives an AFLBR. Indeed, $P^{n_k}(f^k \geq -\frac{1}{k}) \geq P^{n_k}(A_k^c) \to 1$ for $k \to \infty$. Moreover, for $k$ sufficiently large, $P^{n_k}(f^k \geq \frac{\alpha}{4}) \geq P^{n_k}(w^k \geq \frac{\alpha}{2}) - P^{n_k}(A_k) \geq \frac{\alpha}{4}$.

(d) ($\Rightarrow$) Suppose there is an AA2. This implies that there exist $\alpha > 0$ and $f^k \in \mathbf{C}^{n_k}$ with $f^k \geq -1$ and $\lim_{k \to \infty} P^{n_k}(A_k) = 1$, where $A_k = \{f^k \geq \alpha\}$. Define $w^k = \alpha \mathbb{1}_{A_k}$. Then $w^k \in D^{\alpha/2, n_k}$ for all $k$ sufficiently large. For all $(R^n) \triangleleft \triangleright (P^n)$ and for all $u \in \mathbb{U}_4$,

$$E_{R^{n_k}}[u(f^k - w^k)] = E_{R^{n_k}}[u((f^k - w^k) \mathbb{1}_{A_k})] + E_{R^{n_k}}[u((f^k - w^k) \mathbb{1}_{A_k^c})]$$
$$\geq u(-1) R^{n_k}(A_k^c) + u(0) R^{n_k}(A_k).$$

Passing to the limit for $k \to \infty$, we obtain a contradiction to NAMFL($\mathbb{U}_4$) since $R^{n_k}(A_k) \to 1$ for all $(R^n) \triangleleft \triangleright (P^n)$.

(e) Consider the following easy example. Let $1 > \alpha > 0$ and for each $n$, define a random variable $w^n$ on a suitable base $(\Omega^n, \mathcal{F}^n, P^n)$ by

$$w^n = \begin{cases} 1, & \text{on } A^n, \text{ where } P^n(A^n) = \alpha, \\ 0, & \text{on } B^n, \text{ where } P^n(B^n) = 1 - \alpha. \end{cases}$$

The process $(S_t^n)_{t \in \{0,1\}}$ is defined as $S_0^n = 0$ and $S_1^n = w^n$, $\mathcal{F}_0^n$ is trivial and $\mathcal{F}_1^n = \sigma(w^n)$. It is obvious that this large financial market satisfies NAA2 and NSAA, as there cannot exist a subsequence and $\beta > 0$ such that $\lim P^n(w^n \geq \beta) = 1$. On the other hand, there is an AMFL($\mathbb{U}_5$) [and hence an AMFL($\mathbb{U}$) for every subclass $\mathbb{U} \subseteq \mathbb{U}_5$]. Indeed, obviously, $w^n = (S_1^n - S_0^n) \in \mathbf{C}^n$ and moreover, $w^n \in D^{\alpha,n}$. Hence, for every $u \in \mathbb{U}_5$ and every $(R^n) \triangleleft \triangleright (P^n)$,

$$\sup_{f^n \in \mathbf{C}^n} E_{R^n}[u(f^n - w^n)] \geq E_{R^n}[u(w^n - w^n)] = u(0). \qquad \square$$

## APPENDIX

### A.1. Orlicz spaces and Young functions.
Here we follow Kusuoka [19].

DEFINITION A.1.   $F : [0, \infty) \to [0, \infty)$ is a Young function if:

(i) $F$ is continuously differentiable,
(ii) $F(0) = F'(0) = 0$,
(iii) $F'$ is strictly increasing and $\lim_{t \uparrow \infty} F'(t) = \infty$.

Note that the definition of a Young function is not always the same as that as given in Definition A.1—the differentiability is not required in general; see [20]. However, as we will closely follow the approach of Kusuoka [19], we will use his definition of a Young function which is given above. Moreover, it is convenient to use the differentiability in some places.

The class of all Young functions is denoted by $\mathcal{Y}$. For $F \in \mathcal{Y}$, there is a complementary Young function $G \in \mathcal{Y}$ defined as follows:

$$G(y) = \int_0^y (F')^{-1}(t)\,dt, \qquad y \geq 0.$$

Here $(F')^{-1}$ is the inverse function of $F'$. By definition, it is clear that the complementary Young function of the complementary Young function $G$ of $F$ is again the original $F$. The complementary Young function satisfies

(A.1)                          $$G(y) = \max_{x \geq 0}(xy - F(x)).$$

For each $F \in \mathcal{Y}$, let

$$L_F(P) = \{f \in L^0(\Omega, \mathcal{F}, P) : E_P[F(a|f|)] < \infty \text{ for some } a > 0\}.$$

Define a norm $\|\cdot\|_F$ on $L_F$ by $\|f\|_F = \inf\{a > 0 : E[F(\frac{1}{a}|f|)] \leq 1\}$. The space $L_F$ is called an *Orlicz space*. Let

$$L_F^0(P) = \{f \in L^0(\Omega, \mathcal{F}, P) : E_P[F(a|f|)] < \infty \text{ for all } a > 0\}.$$

**Lemma A.2.** *Let $F \in \mathcal{Y}$ and $G$ be its complementary Young function. Then:*

(i) *$L_F$ is a Banach space with norm $\|\cdot\|_F$;*

(ii) *if $g \in L_G$, then $\Phi : L_F^0 \to \mathbb{R}$ given by $\Phi(f) = E[fg]$ is a continuous linear functional and*

$$\|g\|_G \leq \|\Phi\|_{(L_F^0)^*} \leq 2\|g\|_G;$$

(iii) *If $f \in L_F$ and $g \in L_G$, the following relation holds:*

$$E[|fg|] \leq 2\|f\|_F \|g\|_G.$$

**A.2. Some properties of Young functions and the connection to $\mathbb{U}_2$.**

**Lemma A.3.** *Let $F \in \mathcal{Y}$. Then the function $\frac{F(x)}{x}$ is strictly increasing. Moreover, $\lim_{x\to 0} \frac{F(x)}{x} = 0$ and $\lim_{x\to\infty} \frac{F(x)}{x} = +\infty$.*

**Proof.** Let $0 \leq x_1 < x_2$ and define $\lambda = \frac{x_1}{x_2} < 1$. The strict convexity of $F$ and the fact that $F(0) = 0$ imply that $F(\lambda x_2) < \lambda F(x_2)$, which shows the monotonicity of $\frac{F(x)}{x}$. Moreover, by l'Hospital's rule,

$$\lim_{x\to 0} \frac{F(x)}{x} = \lim_{x\to 0} \frac{F'(x)}{1} = F'(0) = 0 \quad \text{and} \quad \lim_{x\to\infty} \frac{F(x)}{x} = \lim_{x\to\infty} \frac{F'(x)}{1} = \infty.$$

$\square$

Let $F \in \mathcal{Y}$. Define a function $u_F \in \mathbb{U}_2$ as follows:

$$\text{(A.2)} \qquad u_F(x) = \begin{cases} -F(-x), & \text{for } x \leq 0, \\ 0, & \text{for } x > 0. \end{cases}$$

**Lemma A.4.** *Let $u \in \mathbb{U}_2$ with $u(0) = 0$. For each $\varepsilon > 0$, there exists $F^\varepsilon \in \mathcal{Y}$ such that $u_{F^\varepsilon}(x) - \varepsilon \leq u(x)$ for all $x \in \mathbb{R}$, where $u_{F^\varepsilon}$ is defined as in (A.2).*

**Proof.** The only problem occurs if the left derivative of $u$ at $0$ is strictly positive. Let $\delta$ be such that $u(-\delta) \geq -\varepsilon$. Define $u_\varepsilon$ as above and such that $u_\varepsilon(x) \leq u(x)$ for $x < -\delta$ (it is clear that there exists $F_\varepsilon \in \mathcal{Y}$ such that this works). Moreover, one can choose $u_\varepsilon$ such that $0 \geq u_\varepsilon(x) \geq u(x)$ on the interval $[-\delta, 0]$. So, $|u_\varepsilon(x) - u(x)| \leq |u(-\delta)| \leq \varepsilon$ for $x \in [-\delta, 0]$. It is clear that $u_\varepsilon \in \mathbb{U}_2$. $\square$

Lemma A.5.   *If $u_F$ is defined as in (A.2), then the convex dual*

$$v(y) = \sup_{x \in \mathbb{R}} (u_F(x) - xy), \qquad y \geq 0,$$

*is in $\mathcal{Y}$. In particular, $v$ is the complementary Young function of $F$.*

Proof.   By the definition of $u_F$,

$$v(y) = \sup\left(\sup_{x>0}(-xy), \sup_{x \leq 0}(-F(-x) - x)\right)$$

$$= \sup_{x \leq 0}(-F(-x) - xy)$$

$$= \sup_{x \geq 0}(-F(x) + xy).$$

This is the definition of the complementary Young function; see (A.1).   □

**A.3. Orlicz spaces and the Mackey topology.**   Define the closed unit ball $B^F(P)$ of the Orlicz space $L_F(P)$ as follows:

$$(A.3) \qquad\qquad B^F(P) = \{f \in L_F(P) : \|f\|_F \leq 1\}.$$

Lemma A.6.   *The set $B^F(P)$ is closed in $L^1(P)$ and uniformly integrable. In particular, $\sup_{h \in B^F(P)} E[|h|\mathbb{1}_{\{|h| \geq \kappa\}}] \leq \frac{\kappa}{F(\kappa)}$.*

Proof.   For $h \in B^F$ and $\kappa > 0$,

$$E[|h|\mathbb{1}_{\{|h| \geq \kappa\}}] = E\left[\frac{|h|}{F(|h|)} F(|h|)\mathbb{1}_{\{|h| \geq \kappa\}}\right] \leq \frac{\kappa}{F(\kappa)} E[F(|h|)] \leq \frac{\kappa}{F(\kappa)},$$

since by Lemma A.3, $\frac{F(y)}{y}$ is strictly increasing. As $\lim_{y \to \infty} \frac{F(y)}{y} = +\infty$,

$$\lim_{\kappa \uparrow \infty} \sup_{h \in B^F(P)} E_P[|h|\mathbb{1}_{\{|h| \geq \kappa\}}] \leq \lim_{\kappa \uparrow \infty} \frac{\kappa}{F(\kappa)} = 0,$$

and this shows uniform integrability.

To show closedness in $L^1$, take $h^n \in B^F$ with $h^n \to h$ in $L^1$. Then a subsequence of $F(|h^n|)$ (still denoted by $n$) converges to $F(|h|)$ a.s. By Fatou's lemma,

$$E[F(|h|)] = E[\lim F(|h^n|)] \leq \liminf E[F(|h^n|)] \leq 1$$

since $F(|h^n|) \geq 0$. This shows that $h \in B^F$.   □

Moreover, the following result holds; see [19].

Lemma A.7.   *Let $A$ be a uniformly integrable subset of $L^1$. Then there exists $F \in \mathcal{Y}$ such that $A \subseteq B^F$.*

So, the uniformly integrable subsets of $L^1$ can be described completely with the help of the sets $B^F$, where $F \in \mathcal{Y}$. Define now the polar of $B^F$:

(A.4)  $V^F(P) = (B^F)^\circ = \{g \in L^\infty(P) : |E_P[gh]| \leq 1 \text{ for all } h \in B^F(P)\}$.

REMARK A.8.  Note that $B^F$ is a weakly compact, convex, balanced subset of $L^1$. Therefore, $B^F = (B^F)^{\circ\circ} = (V^F)^\circ$ by the bipolar theorem; see [7].

Indeed, convexity is clear and so $L^1$-closedness implies closedness for the topology $\sigma(L^1(P), L^\infty(P))$ (i.e., the weak topology on $L^1$). Together with uniform integrability, this gives compactness for the weak topology; see for example [23]. It is clear that $B^F$ is balanced.

The following description of the Mackey topology on $L^\infty$ is used in Section 5. The sets $V^F$, for all $F \in \mathcal{Y}$, give a fundamental system for all Mackey neighborhoods of 0. Indeed, $V^F$ is the polar of the weakly compact subset $B^F(P)$ of $L^1(P)$. On the other hand, by Lemma A.7, it is clear that for each uniformly integrable set in $L^1(P)$, there exists $G \in \mathcal{Y}$ such that $U \subseteq B^G$. So, the polars of all $B^F(P)$ form a fundamental system of neighborhoods for the Mackey topology [which is the topology of uniform convergence on all weakly compact sets of $L^1(P)$]. For the details on the Mackey topology, see [7]. The following lemma gives a helpful description of the sets $V^F$:

LEMMA A.9.  *Let $G$ be the complementary Young function of $F \in \mathcal{Y}$. Then $\frac{1}{2} B^G \cap L^\infty(P) \subseteq V^F \subseteq B^G \cap L^\infty(P)$.*

PROOF.  As, by Lemma A.2, $E_P[|hg|] \leq 2\|h\|_F \|g\|_G$, the first inequality is clear. As for the second inequality, let $g \in V^F$. As $g \in L^\infty(P)$, it is clear that $g$ defines a continuous linear functional $\Phi$ on $L^0_F(P)$ via $\Phi(\cdot) = E_P[g\cdot]$. As $|E[gh]| \leq 1$ for all $h \in B^F(P)$, it is clear that $\|\Phi\|_{(L^0_F)^*} \leq 1$. Therefore $\|g\|_G \leq \|\Phi\|_{(L^0_F)^*} \leq 1$, again by Lemma A.2. This gives that $g \in B^G \cap L^\infty(P)$.  □

**A.4. Contiguity of sequences of probability measures.**  Let $P^n$ and $Q^n$ be probability measures on a probability space $(\Omega^n, \mathcal{F}^n)$. The concept of contiguity can be interpreted as uniform absolute continuity. For a precise statement, see Lemma A.10 below; the easy proof can be found in, for example, [12].

LEMMA A.10.  *Let $Q^n \ll P^n$ for all $n$. Then $(Q^n) \triangleleft (P^n)$ is equivalent to the following condition. For all $\varepsilon > 0$, there exists $\delta > 0$ such that for all $n \in \mathbb{N}$, $Q^n(A^n) < \varepsilon$ for all $A^n \in \mathcal{F}^n$ with $P^n(A^n) < \delta$.*

Moreover, contiguity can be related to uniform integrability. The following result can be found in a more general form in [9]:

LEMMA A.11. *Assume that for all $n \in \mathbb{N}$, $Q^n \ll P^n$. Then $(Q^n) \triangleleft (P^n)$ if and only if $(\frac{dQ^n}{dP^n} | P^n)$ is uniformly integrable, that is,*

$$\lim_{\kappa \uparrow \infty} \sup_n \mathbb{E}_{P^n}\left[\frac{dQ^n}{dP^n}\mathbb{1}_{\{dQ^n/dP^n > \kappa\}}\right] = 0.$$

Together with Lemma A.7, this gives the following characterization of contiguity in terms of Young functions/Orlicz spaces:

LEMMA A.12. *Let $Q^n \ll P^n$ for all $n$. Then $(Q^n) \triangleleft (P^n)$ if and only if there exists $F \in \mathcal{Y}$ such that $\frac{dQ^n}{dP^n} \in B^F(P^n)$ for all $n$.*

PROOF. The sets $B^F(P^n)$ satisfy the uniform integrability condition

$$\lim_{\kappa \uparrow \infty} \sup_n \sup_{h^n \in B^F(P^n)} E_{R^n}[h^n \mathbb{1}_{\{h^n \geq \kappa\}}] = 0;$$

see Lemma A.6. The rest is clear by Lemma A.11 above and Lemma A.7. □

LEMMA A.13. *Suppose that for each $\varepsilon > 0$, there exists a sequence $(Q^{n,\varepsilon})$ with $Q^{n,\varepsilon} \in M^n$ such that:*

(i) *there exists $\delta > 0$ such that for $A^n \in \mathcal{F}^n$ with $P^n(A^n) \geq \varepsilon$, it holds that $Q^{n,\varepsilon}(A^n) \geq \delta$ for all $n$ and*
(ii) *$(Q^{n,\varepsilon}) \triangleleft (P^n)$.*

*Let $Q^n = \sum_{j=1}^{\infty} 2^{-j} Q^{n,2^{-j}}$. Then $Q^n \in M^n \cap \mathbb{P}^n$ for all $n$ and, moreover, $(Q^n) \triangleleft \triangleright (P^n)$.*

PROOF. $\{\frac{dQ^n}{dP^n} = 0\} \subseteq \bigcap_{j=1}^{\infty} \{\frac{dQ^{n,2^{-j}}}{dP^n} < \delta_j\}$, so (i) implies that for all $n$, $P^n[\frac{dQ^n}{dP^n} = 0] = 0$. Hence, $Q^n \in M^n \cap \mathbb{P}^n$.

Let us now prove that $(Q^n)_{n \geq 1} \triangleleft (P^n)_{n \geq 1}$. By Lemma A.10, one must show that for all $\gamma > 0$, there exists $\mu > 0$ such that for all $n$, $Q^n(A^n) < \gamma$ for all $A^n \in \mathcal{F}^n$ with $P^n(A^n) < \mu$.

Let $\gamma > 0$ be fixed and $N \in \mathbb{N}$ sufficiently large so that $\sum_{j=N+1}^{\infty} 2^{-j} < \frac{\gamma}{2}$. By (ii), $(Q^{n,2^{-j}})_{n \geq 1} \triangleleft (P^n)_{n \geq 1}$ for $j = 1, 2, \ldots, N$, whence for each $j$, there exists $\mu_j > 0$ such that for all $n$, whenever $P^n(A^n) < \mu_j$ it follows that $Q^{n,2^{-j}}(A^n) < \frac{\gamma}{2}$. Now let $\mu = \min_{j \leq N} \mu_j$ and $A^n \in \mathcal{F}^n$ be such that $P^n(A^n) < \mu$. Then

$$Q^n(A^n) = \sum_{j=1}^N 2^{-j} Q^{n,2^{-j}}(A^n) + \sum_{j=N+1}^{\infty} 2^{-j} Q^{n,2^{-j}}(A^n) < \frac{\gamma}{2} + \frac{\gamma}{2} = \gamma.$$

To prove that $(P^n)_{n \geq 1} \triangleleft (Q^n)_{n \geq 1}$, observe that (i) implies that for all $j \in \mathbb{N}$, there exists $\mu_j$ such that for all $n$, $Q^{n,2^{-j}}(A^n) < \mu_j$ implies that $P^n(A^n) < 2^{-j}$. Let $\gamma > 0$ be fixed and choose $N \in \mathbb{N}$ such that $2^{-(N-1)} < \gamma$. Define $\mu = 2^{-2N}\mu_N$. Now let $A^n \in \mathcal{F}^n$ be such that $Q^n(A^n) < \mu$. Then

$$P^n[A^n] = P^n\left[A^n \cap \left\{\frac{dQ^{n,2^{-N}}}{dP^n} < \mu_N\right\}\right] + P^n\left[A^n \cap \left\{\frac{dQ^{n,2^{-N}}}{dP^n} \geq \mu_N\right\}\right]$$

$$< 2^{-N} + P^n\left[A^n \cap \left\{\frac{dQ^n}{dP^n} \geq 2^{-N}\mu_N\right\}\right] < 2^{-N} + \frac{2^N}{\mu_N}Q^n(A^n) < \gamma.$$

□

**A.5. Two useful results.** The following result of Bellini and Frittelli (Theorem 2.1 and Corollary 2.2 in [1]) is based on [21]:

THEOREM A.14. $u : \mathbb{R} \to \mathbb{R}$ is nondecreasing and concave on $\mathbb{R}$ (i.e., $u \in \mathbb{U}_2$). Let $R \sim P$ and $G$ be a convex cone with $-L_+^\infty \subseteq G \subseteq L^\infty$. Let $w \in L^\infty$ be such that

$$\sup\{E_R[u(f + w)] : f \in G\} < u(+\infty).$$

Then $N = \{Q \ll R : \frac{dQ}{dR} \in G^\circ \cap L^1(R)\} \neq \varnothing$ and

$$\sup\{E_R[u(f + w)] : f \in G\} = \min_{Q \in N} \min_{\lambda > 0} E_R\left[\lambda \frac{dQ}{dR}w + v\left(\lambda \frac{dQ}{dR}\right)\right],$$

where $v(y) = \sup_{x \in \mathbb{R}}(u(x) - xy)$, $G^\circ = \{\xi \in ba(\Omega, \mathcal{F}, R) : \xi(f) \leq 0 \ \forall f \in G\}$ is the polar of $G$ and $ba(\Omega, \mathcal{F}, R)$ is the space of all bounded, finitely-additive measures on $(\Omega, \mathcal{F})$ which are absolutely continuous with respect to $R$.

In Section 4, Proposition 1.3 of [17] is applied. For the convenience of the reader it is recalled here:

PROPOSITION A.15 (Quantitative version of the Halmos–Savage theorem). For fixed $\varepsilon > 0$ and $\delta > 0$, the following statement is true. Let $M$ be a convex set of $P$-absolutely continuous probability measures such that for all sets $A \in \mathcal{F}$ with $P(A) > \varepsilon$, there exists $Q \in M$ with $Q(A) > \delta$. Then there exists $Q_0 \in M$ such that $Q_0(A) > \frac{\varepsilon^2 \delta}{2}$ for all $A \in \mathcal{F}$ with $P(A) > 4\varepsilon$.

## REFERENCES

[1] BELLINI, F. and FRITTELLI, M. (2002). On the existence of minimax martingale measures. *Math. Finance* **12** 1–21. MR1883783

[2] DALANG, R. C., MORTON, A. and WILLINGER, W. (1990). Equivalent martingale measures and no arbitrage in stochastic securities market models. *Stochastics Stochastics Rep.* **29** 185–201. MR1041035

[3]  DELBAEN, F. and SCHACHERMAYER, W. (1994). A general version of the fundamental
     theorem of asset pricing. *Math. Ann.* **300** 463–520. MR1304434

[4]  DELBAEN, F. and SCHACHERMAYER, W. (1998). The fundamental theorem of
     asset pricing for unbounded stochastic processes. *Math. Ann.* **312** 215–250.
     MR1671792

[5]  FRITTELLI, M. (2004). Some remarks on arbitrage and preferences in securities mar-
     ket models. *Math. Finance* **14** 351–357. MR2070168

[6]  HARRISON, J. M. and PLISKA, S. R. (1981). Martingales and stochastic integrals
     in the theory of continuous trading. *Stochastic Process. Appl.* **11** 215–260.
     MR0622165

[7]  HORVATH, J. (1966). *Topological Vector Spaces and Distributions.* Addison–Wesley,
     Reading, MA. MR0205028

[8]  HUBERMAN, G. (1982). A simple approach to arbitrage pricing theory. *J. Econom.
     Theory* **28** 183–191.

[9]  JACOD, J. and SHIRYAEV, A. N. (1980). *Limit Theorems for Stochastic Processes.*
     Springer, Berlin. MR0959133

[10] KABANOV, Y. and KRAMKOV, D. O. (1994). Large financial markets: Asymptotic
     arbitrage and contiguity. *Theory Probab. Appl.* **39** 222–229. MR1348197

[11] KABANOV, Y. and KRAMKOV, D. O. (1998). Asymptotic arbitrage in large financial
     markets. *Finance Stoch.* **2** 143–172. MR1806101

[12] KLEIN, I. (2000). A fundamental theorem of asset pricing for large financial markets.
     *Math. Finance* **10** 443–458. MR1785165

[13] KLEIN, I. (2003). Free lunch for large financial markets with continuous price pro-
     cesses. *Ann. Appl. Probab.* **13** 1494–1503. MR2023885

[14] KLEIN, I. (2005). No asymptotic free lunch reviewed in the light of Orlicz spaces.

[15] KLEIN, I. (2006). A comment on market free lunch and free lunch. *Math. Finance* **16**
     583–588. MR2239593

[16] KLEIN, I. and SCHACHERMAYER, W. (1996). Asymptotic arbitrage in non-complete
     large financial markets. *Theory Probab. Appl.* **41** 927–934. MR1687136

[17] KLEIN, I. and SCHACHERMAYER, W. (1996). A quantitative and a dual version of
     the Halmos–Savage theorem with applications to mathematical finance. *Ann.
     Probab.* **24** 867–881. MR1404532

[18] KREPS, D. M. (1981). Arbitrage and equilibrium in economies with infinitely many
     commodities. *J. Math. Econom.* **8** 15–35. MR0611252

[19] KUSUOKA, S. (1993). A remark on arbitrage and martingale measure. *Publ. Res. Inst.
     Math. Sci.* **29** 833–840. MR1245018

[20] RAO, M. M. and REN, Z. D. (1991). *Theory of Orlicz Spaces.* Dekker, New York.
     MR1113700

[21] ROCKEFELLAR, R. T. (1971). Integrals which are convex functionals. II. *Pacific J.
     Math.* **39** 439–469. MR0310612

[22] ROSS, S. A. (1976). The arbitrage theory of asset pricing. *J. Econom. Theory* **13**
     341–360. MR0429063

[23] WOJTASZCZYK, P. (1991). *Banach Spaces for Analysts.* Cambridge Univ. Press.
     MR1144277

Institute of Statistics
  and Decision Support Systems
Brünnerstrasse 72
A-1210 Wien
Austria
E-mail: Irene.Klein@univie.ac.at